\documentclass{amsproc}
\usepackage{amsmath,amssymb,amsfonts}
\usepackage{eucal}

\newtheorem{theorem}{Theorem}
\newtheorem{lemma}[theorem]{Lemma}
\newtheorem{prop}[theorem]{Proposition}

\theoremstyle{definition}
\newtheorem{definition}[theorem]{Definition}

\theoremstyle{remark}

\numberwithin{equation}{section}

\input xy
\xyoption{all}

\newcommand{\Map}{\text{Map}}
\newcommand{\SSets}{\mathcal{SS}ets}
\newcommand{\SSetst}{\mathcal{SS}ets^{\mathcal T}}
\newcommand{\SSetsd}{\mathcal{SS}ets^{\mathcal D}}

\newcommand{\LSSetst}{\mathcal {LSS}ets^{\mathcal T}}
\newcommand{\LSSetstm}{\mathcal {LSS}ets^{\mathcal T_M}}
\newcommand{\Tm}{\mathcal T_M}

\newcommand{\Algt}{\mathcal Alg^\mathcal T}
\newcommand{\Algtm}{\mathcal Alg^{\mathcal T_M}}

\begin{document}

\title{Correction to ``Simplicial monoids and Segal categories"}

\author{Julia E. Bergner}
\address{Kansas State University, 138 Cardwell Hall, Manhattan, KS 66506}
\email{bergnerj@member.ams.org}

\maketitle

In the paper \cite{simpmon}, there is an error in the statement of
Proposition 3.12.  The given pair of maps is not in fact an
adjoint pair, since the left adjoint does not preserve coproducts.
Therefore, we give the following revised statement.

\begin{prop}
There is a Quillen equivalence of model categories
\[ \xymatrix@1{L:\LSSetstm_* \ar@<.5ex>[r] & \Algtm:N. \ar@<.5ex>[l]} \]
\end{prop}

The right adjoint functor $N$ is just the forgetful functor; if we
regard a $\Tm$-algebra as strictly local $\Tm$-diagram of
simplicial sets, then this functor forgets the strictly local
structure.

This proposition is then just a variation on Badzioch's
rigidification theorem \cite{bad}.

We give a proof of this proposition, not just for the theory $\Tm$
of monoids, but for any (possibly multi-sorted) algebraic theory,
following the proof of the generalization of Badzioch's
rigidification theorem found in \cite{multisort}.

\begin{prop} \label{Prop}
Let $\mathcal T$ be a (multi-sorted) algebraic theory.  Then there
is a Quillen equivalence of model categories between $\Algt$ and
$\LSSetst_\mathcal O$.
\end{prop}

The proof of the result that we would like to use is given by
letting $\mathcal T =\Tm$, the theory of monoids, and letting
$\mathcal O = \ast$, the set with a single element.

Here, if $\mathcal T$ is an $\mathcal O$-sorted theory,
$\LSSetst_\mathcal O$ denotes the category of functors $\mathcal T
\rightarrow \SSets$ with the homotopy $\mathcal T$-algebra model
structure, but with the additional condition that the image of the
terminal object in $\mathcal T$ is actually isomorphic to the
constant simplicial set given by $\mathcal O$, rather than just
weakly equivalent to it.

We need to find an adjoint pair of functors between $\Algt$ and
$\LSSetst_\mathcal O$ and prove that it is a Quillen equivalence.
Let
\[ J_\mathcal T: \Algt \rightarrow
\mathcal{SS}ets^\mathcal T_\mathcal O \] be the inclusion functor.
We need to show we have an adjoint functor taking an arbitrary
diagram in $\SSetst$ to a $\mathcal T$-algebra.  Here, we use the
idea that a $\mathcal T$-algebra is a strictly local diagram, as
given by the following definition.

\begin{definition}
Let $\mathcal D$ be a small category and $\SSetsd$ the category of
functors $\mathcal D \rightarrow \SSets$. Let $P$ be a set of
morphisms in $\SSetsd$.  An object $Y$ in $\SSetsd$ is
\emph{strictly} $P$-\emph{local} if for every morphism $f:A
\rightarrow B$ in $P$, the induced map on function complexes
\[ f^*: \Map (B,Y) \rightarrow \Map (A,Y) \]
is an isomorphism of simplicial sets. A map $g:C \rightarrow D$ in
$\SSetsd$ is a \emph{strict} $P$-\emph{local equivalence} if for
every strictly $P$-local object $Y$ in $\SSetsd$, the induced map
\[ g^*:\Map (D,Y) \rightarrow \Map(C,Y) \]
is an isomorphism of simplicial sets.
\end{definition}

Now, given a category of $\mathcal D$-diagrams in $\SSets$ and the
full subcategory of strictly $P$-local diagrams for some set $P$
of maps, we have the following result which can be proved just as
in \cite[5.6]{multisort}.

\begin{lemma} \label{adjoint}
Consider two categories, the category of strictly local diagrams
with respect to the set of maps $P= \{f:A \rightarrow B\}$, and
the the category of diagrams $X: \mathcal D \rightarrow \SSets$
which are strictly local with respect to only one of the maps in
$P$. Then the forgetful functor from the first category to the
second has a left adjoint.
\end{lemma}

Now, we can apply this lemma using the fact that a strictly local
$\mathcal T$-diagram is precisely a $\mathcal T$-algebra, and
noting that the objects of $\LSSetst_\mathcal O$ are strictly
local with respect to the map specifying the data on the image of
the terminal object of $\mathcal T$.

Applying Lemma \ref{adjoint} to the functor $J_\mathcal T$, we
obtain its left adjoint functor
\[ K_\mathcal T:\SSetst_\mathcal O \rightarrow \Algt. \]  Then the
following proposition follows just as in the general case
\cite[5.9]{multisort}.

\begin{prop}
The adjoint pair of functors
\[ \xymatrix@1{K_\mathcal T: \SSetst_\mathcal O \ar@<.5ex>[r] & \Algt :J_\mathcal T. \ar@<.5ex>[l]} \]
is a Quillen pair.
\end{prop}

Then, we can extend to the localized model structure
$\LSSetst_\mathcal O$ just as in \cite[5.11]{multisort}.

\begin{prop}
The adjoint pair
\[ \xymatrix@1{K_\mathcal T: \LSSetst \ar@<.5ex>[r] & \Algt
:J_\mathcal T \ar@<.5ex>[l]} \] is a Quillen pair.
\end{prop}

The only difference in the proof is that we can remove one map
(the one with respect to which the objects of $\SSetst_\mathcal O$
are already strictly local) as we localize to get
$\LSSetst_\mathcal O$.  With this minor change in the
localization, the proof of Proposition \ref{Prop} follows just as
the proof of the general case \cite[5.13]{multisort}.

Then, the other correction that should be noted is that the model
structure $\LSSetstm$ should be removed from the chain of Quillen
equivalences just preceding Section 5 of \cite{simpmon}.  Although
this chain is not explicitly given for the many-object case in
Section 5, the same changes should be made there.

The author is grateful to the referee of \cite{inv} who pointed
out this error.

\bibliographystyle{amsalpha}

\end{document}